\documentclass{article}
\usepackage{amssymb,amsmath}
\usepackage{graphicx}
\begin{document}

\raggedbottom

\newtheorem{theor}{Theorem} \newtheorem{defin}[theor]{Definition}
\newtheorem{corol}[theor]{Corollary} \newtheorem{prop}{Proposition}
\newtheorem{lem}[theor]{Lemma} \newcommand{\blem}{\begin{lem}}
\newcommand{\elem}{\end{lem}} \newcommand{\bp}{\begin{prop}}
\newcommand{\ep}{\end{prop}} \newtheorem{example}{Example}
\newcommand{\bex}{\begin{example}} \newcommand{\eex}{\end{example}}
\newcommand{\sq}{\lhd\!\!\!\rhd} \newcommand{\bt}{\begin{theor}}
\newcommand{\et}{\end{theor}}   \newcommand{\bd}{\begin{defin}}
\newcommand{\ed}{\end{defin}}  \newcommand{\bco}{\begin{corol}}
\newcommand{\eco}{\end{corol}}   \newcommand {\6}{\\[.6em]}  \newcommand
{\2}{\\[-2em]}  \newcommand {\be}{\begin{enumerate}}  \newcommand
{\ee}{\end{enumerate}}    
 \newcommand {\bi}{\begin{itemize}}  \newcommand{\ei}{\end{itemize}}   
 \newcommand{\I}{\item}

\begin{center}

\noindent 
{\huge{\bf Defining Homomorphisms and Other Generalized Morphisms of Fuzzy
 Relations in Monoidal Fuzzy Logics by Means of BK-Products.} } 
%\\[1em]} 
\vspace*{0.2in}

\textbf{{\fontsize{11}{12}\selectfont Ladislav J. Kohout}}, \\
\textit{Dept. of Computer Science, Florida State University,
Tallahassee, Florida 32306-4530, USA.\\
E-mail: kohout@cs.fsu.edu \\
%URL: http://cs.fsu.edu/$\!\sim \!$~kohout 
}
\end{center}
%\begin{multicols}{2}
\begin{abstract}
We generalize the previous results that were obtained by Kohout
for relations based on fuzzy Basic Logic systems (BL) of H{\'{a}}jek
and also for relational systems based on left-continuous t-norms. The
present paper  extends  generalized morphisms  into the realm of
Monoidal Fuzzy Logics by first proving and  then using relational
inequalities over pseudo-associative BK-products of relations in these
logics.  \\[1em]
\textbf{Keywords:} \textit{
BK-products of relations, Generalized morphisms, Fuzzy relations,
Monoidal fuzzy logics, t-norms, MV-algebras, Quantum logics,
Relational inequalities, Residuated lattices, Non-associative
compositions of relations. }\\  
\end{abstract}

\noindent 
\begin{quote}
\textit{Invited and refereed paper presented at 
JCIS 2003 - 7th Joint Conf. on Information Sciences (Subsection: 9th Internat. Conf. on Fuzzy Theory and Technology), Cary, North Carolina, USA; September 2003}\end{quote}

\tableofcontents

\section{Introduction}
Homomorphisms play an important role in mathematics, general system
theory and computing as well as in large number practical
applications that require comparison of structures and their matching.  
Many diverse problems of compatibility of structures can be 
unified by generalizing the concept of a homomorphism. Homomorphisms
have been successfully generalized and form one of the basic concepts 
of mathematics of fuzzy sets. In computing and information sciences we
deal to with heterogeneous relations and one way compatibilities for which
both-ways commutativity of diagrams of mappings are severely
inadequate. That is true for both crisp and fuzzy homomorphisms. 
 In 1977 Bandler and Kohout introduced {\it generalized 
homomorphism, proteromorphism} and {\it amphimorphism}, forward and 
backward compatibility of relations, and non-associative and 
pseudo-associative products (compositions) of relations in crisp
setting \cite{77.rel}. These non-associative products were extended to
fuzzy realm in 1978 \cite{78.5}.  
The proofs in the original papers of Bandler 
and Kohout were based on residuation without specific use of
negation \cite{86.2}. Hence they are generally valid in fuzzy
relational calculi  
based on residuated fuzzy logics.  Hence the concepts of generalized
morphisms, compatibility etc. can be rigorously extended  
to relations based on any system of fuzzy logic using continuous
residuated t-norms.  Rigorous proofs in the first order predicate
calculus for  BL family of fuzzy logics of Petr H{\'{a}}jek
were given by Kohout \cite{98.2}. Kohout \cite{98.8} has also
demonstrated that these relational calculi extend outside BL  to the  systems based on\textit{left-continuous}  t-norm family of fuzzy logics. 

BL systems include the well known G{\"{o}}del, {\L}ukasiewicz and
product systems of  fuzzy logics \cite{hajek.bl}. Algebraicaly,
{\L}ukasievicz system is an instance of MV-algebras, which also have
application in development of quantum logics and measures 
\cite{riecan.integral}. Intuitionistic logics 
and linear logics play a role in theoretical computer science. All
these systems are special instances of monoidal systems of fuzzy
logics pioneered by H{\"{o}}hle. In this paper we provide proofs at
the level of monoidal logics, hence our theory of many-valued logic
based relations subsumes theories of relations in the above quoted
systems. The general picture of the hierarchy of fuzzy logic is
depicted in Figure 1.

\begin{figure}[hh]
%\begin{diagram}[width=6em]
%A &\rTo^{R}& A\\
%\dTo^{F}& & \dTo_{F}\\
%C & \rTo^{S} & C
%\end{diagram}
\center
 \includegraphics[scale=0.7]{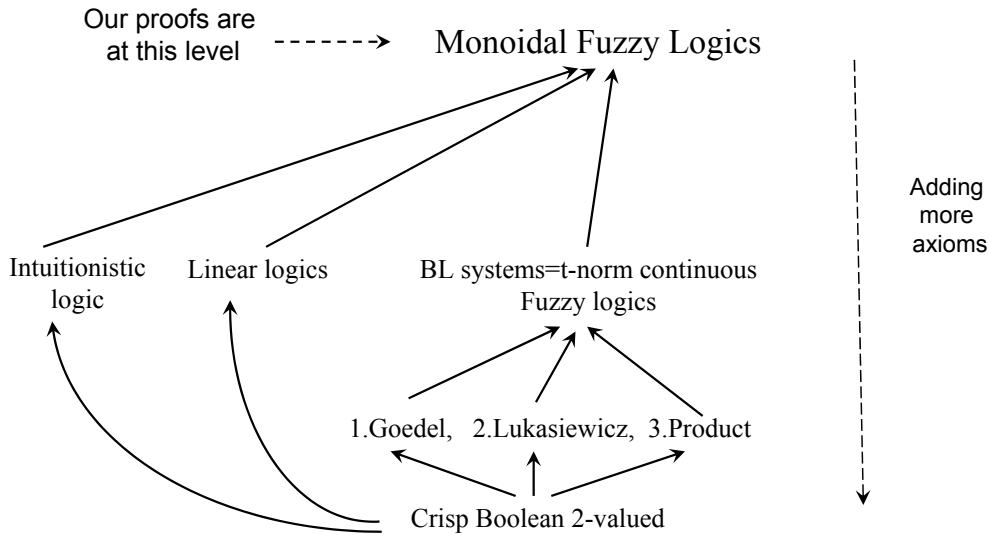}
%\caption{Structure preserving mappings.}
\caption{\textbf{Hierarchy of fuzzy logics of increasing generality}}
\label{fig5}
\end{figure}

\section{Motivation - Crisp Generalized Morphisms} 
\subsection{Crisp Standard Homomorphisms}
Let $A$, $B$, $C$, 
$D$ be sets with relations $R$, $S$ upon them -- $R$ 
from $A$ to $B$ and $S$ from $C$ to $D$, where each relation determines 
some structure. In addition, we have homomorphic mappings $F$ and $G$.  
$F$ is from $A$ to $C$ and $G$ is from $B$ to $D$.  
There are two points of departure that stem from this fundamental 
algebraic notion of {\it homomorphism}: (i) the design or checking 
mappings which will ``preserve" or ``respect" certain given relations, 
and on the other hand (ii) the design or checking of relations which 
``absorb" or ``validate" certain given mappings. 
For example let $A=B$, $C=D$ and $R$, $S$ be orders. Given $A$ and 
$C$ we wish to find one or all the mappings from $A$ to $C$ that 
preserve orders -- this illustrates the case (i). An example of (ii) 
is, given a mapping from $A$ to $C$, how to match the order on $A$ 
given by $R$, with some other order on $C$, or vice versa -- so that 
some given mapping will preserve or co-preserve them. 
Another example is where $A=B \times B$, $C=D \times D$ and $R$ 
and $S$ determine some groupoids 

In this situation, the conventional homomorphism yields a commuting 
diagram of arrows such that $R \circ G = F \circ S$, where of course,  
the morphisms $F$ and $G$ are the relations which are  both covering 
and univalent (i.e. functional). To obtain the constructions 
that solve the problems (i) or (ii) requires to solve the above 
relational equation with respect to one of the relations $R$, $S$, 
$F$ or $G$.  
%FIG1*********

\begin{figure}[hh]
%\begin{diagram}[width=6em]
%A &\rTo^{R}& A\\
%\dTo^{F}& & \dTo_{F}\\
%C & \rTo^{S} & C
%\end{diagram}
\center
 \includegraphics[scale=0.9]{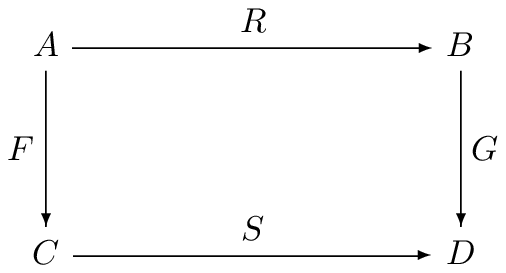}
%\caption{Structure preserving mappings.}
\caption{\it A diagram for conventional homomorphisms.
Here homomorphic maps $F$ and $G$ are functions, i.e. univalent and
covering relations.}
\label{fig1}
\end{figure} 

When the mappings (functional relations) $F$ and $G$ are replaced by 
general relations,  the equation is no longer valid but has to be 
replaced by two inequalities. The notion of a homomorphism splits 
into two independent notions, generalized morphism and generalized 
proteromorphism. 

\subsection{Crisp Generalized Morphisms} 
It is useful to summarize here the basic notions concerning generalized
morphisms of crisp relations, as this information is not generally
available in textbooks despite of the fact that generalized morphisms
and proteromorphisms were introduced by Bandler and Kohout in
1977. Familiarity with the crisp equalities and inequalities
characterizing these will facilitate understanding of the fuzzy
case.

\subsubsection{Partial and Total Homomorphisms}
The following simple observation and Lemma 1 will help to comprehend
the effect  of
relaxing equational constraints defining homomorphisms into
inequalities that characterize generalized morphisms
and generalized proteromorphisms.  

Trivially, 
\[R \circ G = F \circ S\;\;\;\mbox{iff}\;\;\; F \circ S=R\circ G.\]
Composing the left hand side of the above expression with 
the inverse of $F$, the relation $F^{-1}$ applied from the left yields 
$F^{-1}\circ R \circ G =F^{-1}\circ F \circ S=E \circ S=S.$ Hence the 
equation $F^{-1}\circ R \circ G=S$ is equivalent to $R \circ G = F \circ S.$

Similarly, composing the right hand side of the above expression with 
the inverse of $F$, the relation $G^{-1}$ applied from the right yields 
$ F \circ S\circ G^{-1} =R \circ G\circ G^{-1}=R\circ E=R.$ Hence the 
equation  $F \circ S \circ G^{-1}=R$ is equivalent to $F \circ S = R \circ G.$

Hence, the following obvious equivalence holds:
\((F^{-1}\circ R \circ G=S) \equiv (F \circ S \circ G^{-1}=R) \equiv (F
\circ S = R \circ G)\equiv (R \circ G = F \circ S) \)\\[1em]  
The diagram of the Fig. 1  is a partial homomorphism if the equality
above holds AND the relations $F$ and $G$ are \textit{univalent} (i.e.
partial functions). It is a homomorphism, if in addition both
relations $F$ and $G$ are covering i.e. (total) functions. 

This is summarized in the following obvious lemma: 
\blem [Homomorphism]
For any pair of relations $R$ and $S$, where $S$ is the homomorphic
image of $R$,  the following conditions simultaneously hold:
\be
\I There exist relations $F$ and $G$ such that the equality \\ \((F^{-1}\circ R \circ G=S) \equiv (F \circ S \circ G^{-1}=R) \) holds, and \\[-1.7em]
\I $F$ and $G$ are both univalent and covering relations. 
\ee
\elem 

\blem [Partial Homomorphism]
If the arrows in Fig. 1 commute, i.e. the equality \(R \circ G = F
\circ S\) holds, and $F$ and $G$  are univalent relations, then $S$ is
a partial homomorphic image of $R$ (i.e. partial homomorphism).
\elem

When the relational equality  \((F^{-1}\circ R \circ G=S)\) on the
left, or the relational  equality \(F \circ S \circ G^{-1}=R\) on the
right in expression (1) of 
Lemma 1 is replaced by the relational inclusion $\sqsubseteq$, the
commuting diagram of Fig. 2 splits into two diagrams (see Fig. 3
below) and the notion of
homomorphism has to be replaced by the notion of generalized morphisms
as described in the next section.  

\subsubsection{Crisp Generalized Morphisms and Proteromorphisms}

When the homomorphic mappings $F$ are $G$ not functions, the diagram of
Fig. 2 does not commute any more, and the homomorphism does not exist in
general case.  
In that case the equality  \((F^{-1}\circ R \circ G=S)\) on the left, or the
equality on the right in expression (1) of Lemma 1 changes into an
inequality. The notion of homomorphism splits into two different
notions,  Generalized Morphism and Proteromorphisms. The diagrams and
inequalities for these are shown in Figures 3a and 3b.

\begin{figure}[h]
\center
  \includegraphics[scale=0.9]{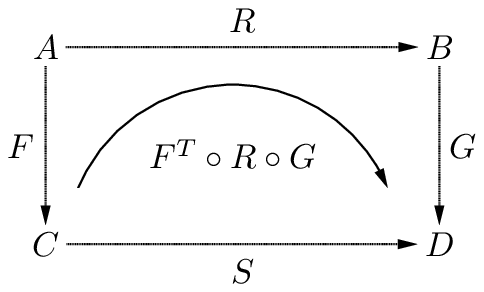}\hspace*{2cm}
 \includegraphics[scale=0.9]{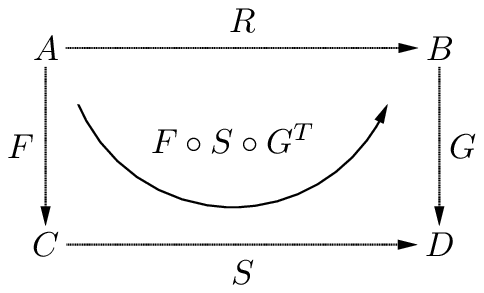} 
\\{\hspace{1.2cm}(a) Generalized homomorphism: \hspace{1.2cm} (b)
 Generalized  Proteromorphism:}
\\{$F^{-1}\circ R \circ G\sqsubseteq S$ \hspace{4.2cm}  $F \circ S
  \circ G^{-1}\sqsubseteq R$}
\caption{\it Forward and Backward Compatible Relations}
\label{fig2}
\end{figure}
%\label{fig2}

\section{Solutions of Relational Inequalities Characterizing
  Generalized Morphisms}   

The proofs in the original papers of Bandler 
and Kohout were based on residuation without specific use of
negation [2]. Hence they are generally valid in fuzzy relational calculi 
based on residuated fuzzy logics.  Hence the concepts of generalized
morphisms, compatibility etc. can be rigorously extended  
to relations based on any system of fuzzy logic in which the
implication operator is the residuum of the AND connective. 

In this section we give just a sampler of selected solutions. 
Kohout extended  the previous results of Bandler and Kohout \cite{86.2}
on generalized morphisms to relations based on fuzzy Basic Logic
systems (BL) of H{\'{a}}jek 
and also to  relational systems based on left-continuous t-norms. In
this section we give just a sampler of selected solutions for $R$ and
$S$.  

The solutions for $F$ and $G$ will be presented in the sequel. 
Sections 4 and 5 extend  generalized morphisms  into the realm of
Monoidal Fuzzy Logics by first proving and  then using relational
inequalities over pseudo-associative BK-products of relations in these
logics.

\subsection{From Crisp to Fuzzy Case}
Relational inequalities displayed in Fig. 3 of Sec. 2.2.2 
give a rigorous mathematical definition of generalized morphisms. If
we want to use generalized morphisms either in pure mathematics or in
applications (such as knowledge engineering, scientific computations
etc.) we need some other theorems describing the properties of
generalized morphisms. 

For example, given any three relations chosen from $R, S, F, G$ we may
wish to compute the fourth remaining unknown one.  In order
to do this, we have to possess the solution of inequalities that
allows us to compute the unknown relation for the known
ones. Compatibility criteria provide solution for either $R$ or
$S$. In latter sections we shall also present the solutions for $F$ and $G$. 
 
\subsubsection{Formulation of Compatibility Criteria}
\begin{enumerate}
\item Forward Compatibility 
$F^{-1} \circ R \circ G \sqsubseteq S$ is fulfilled iff
 $R \sqsubseteq F \lhd S \rhd G^{-1}$

\item Backward Compatibility
 $F \circ S \circ G^{-1} \sqsubseteq R$  is fulfilled iff
 $S \sqsubseteq F^{-1} \lhd R \rhd G$
 
\item Bothways Compatibility is characterized by
        \begin{enumerate}
        \item $F \circ S \circ G^{-1} \sqsubseteq R \sqsubseteq F \lhd S \rhd G^{-1}$   \\
                or equivalently by
        \item $F^{-1} \circ R \circ G \sqsubseteq S \sqsubseteq F^{-1} \lhd R \rhd G$
        \end{enumerate}
\end{enumerate}

The solutions involve non-associative compositions of relations called 
BK-products in the literature. We shall work with the sub-product
$\lhd$ and super-product $\rhd$. Before with proceed with further
technicalities of  the proofs, we 
briefly summarize the basic algebraic facts about $\lhd$ and  $\rhd$. 

Algebraic characterization of the interplay of
the triangle sub-product $\lhd$ and the triangle super-product $\rhd$ with
the standard associative $\circ$ (circle) product \textit{forms 
the algebraic core} on which the subsequent proofs are based. 
%This will be followed by the list of solutions. 
%We list all the solutions first, and then proof their validity in the
%monoidal fuzzy logics.  

\subsubsection{Algebraic Properties of BK-products of Relations}

The power of both crisp and fuzzy relational calculi is substantially
enhanced by introducing non-associative compositions of relations in
addition to the well-known standard circle product $\circ$. These
additional relational  
compositions that we  called triangle and square products
\cite{92.7},\cite{98.8},\cite{00.6},\cite{86.5}  were first introduced
by Bandler and Kohout  in 1977  \cite{77.rel},\cite{87.8},\cite{80.3}
and are referred to as the BK-products  in the literature  
\cite{hajek96.bk},\cite{belohl.bk},\cite{deb+kerMP},\cite{deb+ker93pr},\cite{deb+ker94cut},\cite{99.2}.      

The representational and  computational power of BK-products $\rhd$ and
$\lhd$ resides in their algebraic properties. 
The following mixed pseudo-associativities hold for $\lhd$ \  and \
$\rhd$:   \[ Q \lhd (R \rhd S) = (Q \lhd R) \rhd S,\]
 \[Q \lhd (R \lhd S) =(Q \circ R) \lhd S \]
 \[ Q \rhd (R \rhd S) = Q \rhd (R \circ S).\]     
The interplay of $\circ, \lhd, \rhd$ that is afforded by relaxing the
property of full associativity is essential for enriching the
expressive power of the calculus of relations. The mutual interaction
of these three relational compositions plays a crucial role in
defining the key inequalities of relational calculus. 

One such set of inequalities  called  {\it  Residuation bootstrap of
  BK-products} that plays a crucial role in the development of fuzzy
  relational calculi \cite{00.4} will be proved and used extensively in
  this chapter.  It consists of the following  relational
  inequalities that hold  for arbitrary $V \in \cal{B} (A \rightarrow C)$: 
\[R \circ S \sqsubseteq V \;\; \mbox{iff} \;\;  R \sqsubseteq V
\rhd S^{T} \;\; \mbox{iff} \;\; S \sqsubseteq R^{T} \lhd V 
\]

\subsection{ BK-Products of  Relations} 
We shall briefly summarize the basic notions concerning the
non-associative BK-products of relations. This knowledge is essential
for fuller understanding of the proofs of the inequalities
characterizing the mathematical properties of generalized morphisms
hat are presented id latter sections of this paper.

\subsubsection{A Brief Overview of BK-Products}

\noindent
{\bf Mathematical definitions.} \ \  Where $R$ is a relation
from $X$ to $Y$, and $S$ a relation
from $Y$ to $Z$, a {\it product relation} $R @ S$ is a relation
from $X$ to $Z$, determined by $R$ and $S$. There are several
types of product used to produce product-relations \cite{87.8},
\cite{92.7}. Each product type performs a {\bf different logical
action} on the intermediate sets, as {\it each logical type} of
the product enforces a {\it distinct specific meaning} on the
resulting product-relation $R @ S$.
In the following definitions of the products,
$R_{ij}, S_{jk}$ represent the
fuzzy degrees to which the respective statements $x_{i}Ry_{j}$,
$y_{j}Sz_{k}$ are true.

%Table 1 **********
\begin{table}[h]
\center
\caption{$\circ$-product and non-associative  BK-products of
  relations}\vspace{1em}
\renewcommand{\arraystretch}{1.4}
\setlength\tabcolsep{7pt}
\begin{tabular}{@{}l|l|l } %p{1.8cm}}
\hline\noalign{\smallskip}
\textsc{Product Type} &  \textsc{Many-Valued Logic} & 
\textsc{Set-based Definition}\\
\noalign{\smallskip}
\hline
\noalign{\smallskip}
{\bf Circle} product: & $(R\circ S)_{ik}= {\bigvee}_{j} (R_{ij}
\&S_{jk})$ &$ x(R \circ S)z \;\;\;\;\Leftrightarrow 
\;\;\; xR\;\;  \mbox{intersects}\;\; Sz$ \\
{\bf Sub-product} 
& $(R\lhd S)_{ik}= {\bigwedge}_{j} (R_{ij} \rightarrow  S_{jk})$ &
$ x(R \lhd S)z \;\;\;\Leftrightarrow \;\;\; xR 
\stackrel{\sim}{\subseteq} Sz$  \\ 
{\bf Super-product} & $(R
\rhd S)_{ik}= {\bigwedge}_{j} (R_{ij} \leftarrow S_{jk})$
&$ x(R \rhd S)z \;\;\;\Leftrightarrow
\;\;\; xR \stackrel{\sim}{\supseteq} Sz$ \\
{\bf Square} product: & $(R\; \Box\; S)_{ik}= {\bigwedge}_{j} (R_{ij}
\equiv S_{jk})$
& $\; x(R \Box S)z \;\;\;\Leftrightarrow
\;\;\; xR \cong Sz$  \\
\noalign{\smallskip}
\hline
\noalign{\smallskip}
\end{tabular}\\
\label{Tab1}
\end{table}

There are several different notational forms in which
BK-products can be expressed:\\[-2em]
\begin{enumerate} 
\item the notation shown in Table 1 using the concept of fuzzy set inclusion
and  equality  \cite{80.1},\cite{80.3}.
\item many-valued logic(MVL) based notation, which uses the logic
connectives $\bigwedge$, $\&$, $\rightarrow$ or $\equiv$ which is
also displayed in Table 1. 
\item The tensor notation (not needed in this paper).
\item The fuzzy predicate calculus form (see Table 3 in  Sec. 4.1 below). 
\end{enumerate} 
\noindent 
These four different forms of relational compositions are
logically equivalent under some reasonable logic assumptions,
producing the same mathematical results. 
Distinguishing these forms is, however, important when constructing
fast and efficient computational algorithms \cite{92.7}.

The tensor notation in its presentation abstracts from the display of
the type of MVL connectives shown by logic-based notation. It
preserves, on the other hand, the information about the way the
BK-products were composed from their components. This is important when
we want to keep track of the ways in which several distinct, but
logically equivalent streams of relational computation were constructed.

The logical symbols for the logic connectives {\tt AND} $\&$,
both {\it implications} and the {\it equivalence} in the 
formulas shown in Table 1 represent connectives of some many-valued logic,
{\bf chosen} according to the logic properties of the products required. 
An important special case is when the {\tt AND} connective $\&$ is
represented semantically by a t-norm *.  If the logics are residuated, then
the implications are residua of the t-norm, and the equivalence is a
biresiduum of the t-norm.  

The generic formula
\[{\displaystyle (R @ S)_{ik} :=\ \oplus_{j} (R_{ij} \# S_{jk}),}\]
yields two types of fuzzy relational products. We can replace the
outer connective $\bigoplus$ with $\bigwedge$(defined above) or with
$\frac{1}{|J|}\sum$;\\ 

\begin{center}
%\noindent
${\displaystyle (R @ S)_{ik} :=\ \bigwedge_{j} (R_{ij} \# S_{jk})}$:
{\em Harsh} product, \\[0.5em]
${\displaystyle (R @ S)_{ik} :=\ \frac{1}{|J|}\sum_{j} (R_{ij} \# S_{jk})}$:
{\em Mean} product.
\end{center}
\noindent
By choosing appropriate many-valued logic
 operations for  the logic connectives,
the crisp case extends to a wide variety of
many-valued logic based (fuzzy) relational systems \cite{92.7},
\cite{80.2},\cite{86.5},\cite{87.6},\cite{87.13},\cite{92.7}.
While we often used in our applications the classical {\em min} and
 {\em  max} for {\em t-norm} and {\em t-conorm},
respectively, we applied various MVL implication operators for the
computation of BK-products.
The details of choice of the appropriate many-valued connectives are 
discussed in \cite{80.2},\cite{86.5},\cite{87.6},\cite{87.13},\cite{92.7}. 

\section{Residuation Bootstrap of BK-products in
  Monoidal Fuzzy Logics}
Now, we shall look at the ways of generalizing the Residuated
Bootstrap of BK-products to monoidal fuzzy logics.  
It is sufficient to prove that the Residuated Bootstrap of BK-products
holds in fuzzy monoidal logics. The proofs of the inequalities
characterizing generalized morphisms follow then from the bootstrap
inequalities in the same way as in t-norm based fuzzy logics.

\subsection{Residuated Lattices and Monoidal Logics}

BL systems were based on the idea  that many important
theorems of Zadeh's 
logics on $[0, 1]$ would still hold when {\it min} is replaced by any
continuous t-norm $\&_T$ and 
$\rightarrow_3$ by the corresponding residuated $\rightarrow_T$ implication
operator. The logic systems that employ the pair a t-norm and its residuum
$({\&}_{T}, \; \rightarrow_T)$  as {\it and} and {\it implication}
connectives were called {\it Basic Logics} (BL)
\cite{hajek.bl}. One further extension was with left-continuous
t-norms in which our Residuated Bootstrap inequalities also held. 
Our
theorems, however, will be further generalized and shown to hold in
residuated lattices. These lattices form a foundation of fuzzy logics in
monoidal categories \cite{hohle.monclcat}. 
For logics, in order to possess adequate properties,
complete residuated lattices are usually required.

%\subsubsection{What Are Monoidal Fuzzy Logics}

\begin{defin}{\textbf{Residuated Lattice}} (integral, residuated, commutative l-monoid). \\
A residuated lattice ${\cal L}=(L, \leq, \wedge, \vee, \otimes, \rightarrow,
\mbox{\bf 0}, \mbox{\bf 1})$ is a lattice containing the least element {\bf
0} and the largest element {\bf 1} and the additional two 2-argument
operations  $\otimes$ and $\rightarrow$. $\otimes$ is a commutative monoid for which the
{``residuum''}  $\rightarrow$ is
determined by the Galois correspondence given by the formula  $a \otimes b  \leq c
\Longleftrightarrow a \leq b \rightarrow c$. \\[-0.4cm]
\end{defin}

%\subsubsection{Further Properties of Fuzzy Residuated Logics} 

The following formulas that hold in residuated lattices specified by
Def. 3 will be needed in the sequel. We can see that the lattice
semantics can be translated easily into first order logic formulas of
fuzzy monoidal logics as shown in Table 2. 

%Table 2 ****
\begin{table}[h]
%\center 
\caption{Lattice Semantics of the First Order Formula of Fuzzy
  Monoidal Logics} \vspace{-1em} 
\begin{align*} 
\mathbf{Lattice\;\; Semantics} &  & &
\mathbf{\;\;\;\;First\;\;Order\;\;Logic\;\;Formulas}  
\end{align*}
\begin{align} 
(x \otimes y) \Rightarrow z& \; =\; x \Rightarrow (y \Rightarrow z)&  
(\nu \;\&\; \varphi) \rightarrow \psi & \; =\; \nu
\rightarrow (\varphi \rightarrow \psi) \\[1em] 
(x \otimes \bigvee_{i\in I} y_{i})& = \bigvee_{i\in I} (x\otimes y_{i})& 
 \nu \;\&\; (\exists i) \varphi & = (\exists i) (\nu \;\&\;
 \varphi) \\[0.4em] 
(x \Rightarrow \bigwedge_{i\in I} y_{i})& = \bigwedge_{i\in I}
(x\Rightarrow y_{i})&
(\nu \rightarrow (\forall i) \varphi & = (\forall i) 
 (\nu \rightarrow \varphi ) 
\end{align}
\label{Tab2}
\end{table}

%\end{theor}

Table 3  displays the residuated lattice semantics and first
order syntactic formulas of BK-products. This supplements other forms
of BK-product representations that were given in Sec. 3.2.1, in Table
1.  

%Table 3 ***
\begin{table}[h] 
\center 
\textbf{ 
\caption{ $\circ$-product and   BK-products: Semantics and Syntax}}\vspace{1em}
\renewcommand{\arraystretch}{1.4}
\setlength\tabcolsep{7pt}
\begin{tabular}{@{}ll|c|c } %p{1.8cm}}
\hline\noalign{\smallskip}
\textsc{Product Type} & & \textsc{Residuated Lattice} & 
\textsc{First Order Logic}\\
& &  \textsc{Semantics} &  \textsc{Formulas}\\
\noalign{\smallskip}
\hline
\noalign{\smallskip}
{\bf Circle} product: & $(R\circ S)$ & $\bigwedge_{i}\bigwedge_{k} {\bigvee}_{j} (R_{ij}\otimes S_{jk})$
& \( (\forall x) (\forall z) (\exists y) (xRy \; \& \; ySz) \) \\
{\bf Sub-product:} 
& $(R\lhd S)$ &  $\bigwedge_{i}\bigwedge_{k}{\bigwedge}_{j} (R_{ij} \Rightarrow  S_{jk})$ &
 \( (\forall x) (\forall z) (\forall y) (xRy
  \rightarrow ySz) \) \\
{\bf Super-product:} & $(R
\rhd S)$ & $\bigwedge_{i}\bigwedge_{k} {\bigwedge}_{j} (R_{ij}
\Leftarrow S_{jk})$ &\( (\forall x) (\forall z) (\forall y) (xRy
  \leftarrow  ySz) \)\\
{\bf Square} product: & $(R\; \Box\; S)$ & ${\bigwedge_{i}\bigwedge_{k}\bigwedge}_{j} (R_{ij}
\Leftrightarrow S_{jk})$
& \( (\forall x) (\forall z) (\forall y) (xRy
  \equiv  ySz) \) \\
\noalign{\smallskip}
\hline
\noalign{\smallskip}
\end{tabular}\\
\label{Tab3}
\end{table}

\subsection{The Proof of Residuation Bootstrap of BK-Products in
  Monoidal Fuzzy Algebras}
\noindent
\begin{theor} \ {\it Residuation bootstrap of BK-products 
\cite{00.4}.} \\     
For arbitrary $V \in \mathcal{B} (A \rightarrow C)$,
\[T \circ U \sqsubseteq V \;\; \mbox{iff} \;\;  T \sqsubseteq V
\rhd U^{-1} \;\; \mbox{iff} \;\; U \sqsubseteq T^{-1} \lhd V 
\] 
universally holds in residuated lattice (integral, residuated,
commutative l-monoid)   of Def. 3.  
\end{theor}

\noindent {\it Proof:}\\

\begin{tabular}{ll}
 $T \circ U \sqsubseteq V$ & \\
$\bigwedge_{a}\bigwedge_{c}(a(T \circ U)c\Rightarrow aVc)$ &\\
$\bigwedge_{a}  \bigwedge_{c} (\bigvee_{b}( aTb \otimes
 bUc)\Rightarrow aVc )$ & \\ %by $\circ$ definition\\
$ \bigwedge_{a}  \bigwedge_{b} \bigwedge_{c}(bUc \Rightarrow
 \bigwedge_{a} (bT^{-1}a \Rightarrow aVc))$ & \\ %by 4.1.2-(1) \\
$\bigwedge_{a}  \bigwedge_{c} (bUc \Rightarrow b(T^{-1} \rhd V)c)$ & \\
$U \sqsubseteq T^{-1}\rhd V$ & \\[0.4em]% by $\rhd$ definition\\[0.4em]
\end{tabular}

\noindent
Other parts of the formula of Th. 4 can be easily proved in a
similar way.  

\section{Solutions  of Inequalities in Monoidal Fuzzy Logics}

\subsection{Classification of Generalized Morphisms} 

 Generalized Morphisms $F, G$ from relation $R$ to relation $S$ are
 classified in Table 4. \\

%Table 4 *******
\begin{table}[h] 
\center 
\textbf{ 
\caption{An Overview of Generalized Morphisms}}\vspace{1em} 
\noindent
\begin{tabular}{c|c|c}\hline 
& & \\[-0.9em]
 \textbf{Generalized Morphism}& \textbf{Type
  of Compatibility} & \textbf{Relational Definition} \\
 & & \\[-0.9em]\hline  \hline 
& & \\[-0.7em]
 $F, G$ are gen. homomorphisms &   $FRG$: $S$ are \\ 
from $R$ to $S$ & forward-compatible  & $F^{-1} \circ R
  \circ G \sqsubseteq S$ \\
 & & \\[-0.7em]\hline 
& & \\[-0.7em]
 $F, G$ are gen. homomorphisms &  $FRG$:$S$  are & \\
 from $R$ to $S$ & backward-compatible & $F \circ S \circ G^{-1}
 \sqsubseteq R$\\ 
& & \\[-0.7em]\hline 
& & \\[-0.7em]
 $F, G$ are Gen. amphimorphisms &   $FRG$: $S$ are &  $F^{-1} \circ R
  \circ G \sqsubseteq S$\\
 from $R$ to $S$ & bothways-compatible & and \\
& & $F \circ S \circ G^{-1} \sqsubseteq R$ \\[0.4em] \hline 

\end{tabular} 
\label{Tab4}
\end{table}

\bd [Generalized Amphimorphism]
Simultaneous fulfillment of the conditions of backward compatibility
and forward compatiblity will be expressed as both-ways compatibility
and such a morphism will be called Generalized Amphimorphism. See
Fig. 4.
\ed

\begin{figure}[h] 
 \vspace*{2em}
\center
  \includegraphics[scale=0.9]{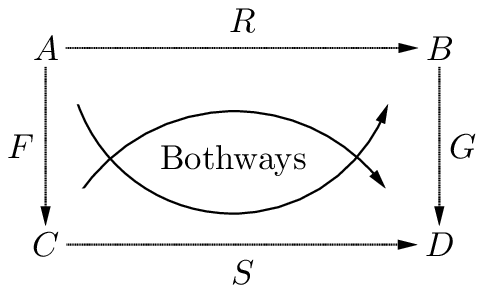}
%\caption{Structure preserving mappings.}
\caption{\it Generalized Ampimorphism or Both-Ways Compatibility. }
\label{fig4}
\end{figure}

Homomorphism is a special kind of both-ways compatibility. 

\subsection{Solutions and  Proofs in Monoidal Fuzzy Algebras}

In the proof we shall use  the Residuation bootstrap of BK-products,
namely 
\[T \circ U \sqsubseteq V \;\; \mbox{iff} \;\;  T \sqsubseteq V
\rhd U^{-1} \;\; \mbox{iff} \;\; U \sqsubseteq T^{-1} \lhd V 
\] 
the validity of which in  Monoidal Fuzzy Algebras we proved in the
previous section (cf. Theorem 4). It will be convenient to split this
expression into two parts denoting these parts as B1 and B2, respectively:

\begin{description}
\I[B1:]\hspace{1.5cm} 
\(T \circ U \sqsubseteq V \;\;\overset{B1}{\Longleftrightarrow}   \;\;
 U \sqsubseteq T^{-1} \lhd V \)  
\I [B2:]\hspace{1.5cm}
\(T \circ U \sqsubseteq V \;\;\overset{B2}{\Longleftrightarrow} \;\;
T \sqsubseteq V \rhd U^{-1}\)
\end{description}

\begin{theor} {\sc Forward Compatibility Solution}. \\
 $FRG:S$ are forward compatible $\Longleftrightarrow$  $F^{T} \circ R 
\circ G \sqsubseteq S$ $\Longleftrightarrow$ 
$R \sqsubseteq F \lhd S \rhd G^{T})$ % \\[-1.4em]
\end{theor}
%\noindent FORWARD:\\ 
%****************************************\\
\textit{Proof}:\\
Substituting \(T:=F^{-1}, \;\; U:=R \circ G, \;\; V:=S \;\;\) into B1
we obtain \\ \( F^{-1}\circ R \circ G \sqsubseteq S \;\;\;
\overset{B1}{\Longleftrightarrow}  
\;\;\;  R \circ G \sqsubseteq F \lhd S \); \\
Substituting \(T:=R, \;\; U:=G, \;\; V:=F \rhd S
\;\;\) into B2 we obtain  \\
\( R \circ G \sqsubseteq F \lhd S \;\;\;
\overset{B2}{\Longleftrightarrow}  \;\;\;  R \sqsubseteq F\lhd S \rhd
G^{-1} \;\;\);   
Transitivity of equivalences yields \(F^{-1}\circ R \circ G \sqsubseteq S
\;\;\;\Longleftrightarrow   
\;\;\;  R \sqsubseteq F \lhd S \rhd G^{-1} \). This completes the
proof.    \\

\begin{theor} {\sc Backward Compatibility Solution}. \\
$FRG:S$ are backward compatible $\Longleftrightarrow$  $F \circ S 
\circ G^{T} \sqsubseteq R$ $\Longleftrightarrow$ 
$S \sqsubseteq F^{T} \lhd R \rhd G$  \\[- 1.6em] 
\end{theor}
%\noindent BACKWARD:\\
\textit{Proof}:\\
Substituting \(T:=F, \;\; U:=S \circ G^{-1}, \;\; V:=R \;\;\) into B1
we obtain \\ \( F\circ S \circ G^{-1} \sqsubseteq R \;\;\;
\overset{B1}{\Longleftrightarrow}  
\;\;\;  S \circ G^{-1} \sqsubseteq F^{-1} \lhd R \); \\
Substituting \(T:=S, \;\; U:=G^{-1}, \;\; V:=F^-{1} \rhd R
\;\;\) into B2 we obtain \\
\( S \circ G^{-1} \sqsubseteq F^{-1} \lhd R \;\;\;
\overset{B2}{\Longleftrightarrow}  \;\;\;  S \sqsubseteq F^{-1} \lhd R \rhd
G \;\;\);  
Transitivity of equivalences yields \( F^{-1}\circ R \circ G \sqsubseteq
S  \Longleftrightarrow  
\;\;\;  S \sqsubseteq F^{-1} \lhd (R \rhd G) \). This completes the
proof.    \\

\begin{theor}[Forward Compatibility: Criteria for $F$ and $G$]
\begin{verbatim}
\end{verbatim}
   $FRG:S$ are forward-compatible iff
        \begin{enumerate}
        \item  $F \sqsubseteq R \lhd (G \lhd S^{-1})$ \\
               or equivalently
        \item  $G \sqsubseteq R^{-1} \lhd (F \lhd S)$ 
        \end{enumerate}
\end{theor}
\textit{Proof:}\\ 
(1): Criterion for F:\\
\( F^{-1}\circ R \circ G \sqsubseteq S \;\;\; \Longleftrightarrow
\;\;\;  F^{-1} \circ R \sqsubseteq S \rhd G^{-1} \;\;\; \Longleftrightarrow
\;\;\;  F^{-1} \sqsubseteq (S \rhd G) \rhd R^{-1} \\ 
 \Longleftrightarrow
\;\;\;  F \sqsubseteq R \lhd (S \rhd G^{-1})^{-1}  \Longleftrightarrow
\;\;\;  F \sqsubseteq R \lhd (G \lhd S^{-1})    \) \\[0.5em]
(2): Criterion for G:\\
\( F^{-1}\circ R \circ G \sqsubseteq S \;\;\; \Longleftrightarrow
\;\;\;  R \circ G \sqsubseteq F \lhd S \;\;\; \Longleftrightarrow
\;\;\;  G \sqsubseteq R^{-1} \lhd (F \rhd S) \)

\begin{theor}[Backward Compatibility: Criteria for $F$ and $G$]
 \begin{verbatim}
\end{verbatim} 
 $FRG:S$ are backward-compatible iff
        \begin{enumerate}
        \item  $F \sqsubseteq (R \rhd G) \rhd S^{-1})$ \\
               or equivalently
        \item  $G \sqsubseteq (R^{-1} \rhd F) \rhd S$
        \end{enumerate}
\end{theor}
\textit{Proof:}\\ 
(1): Criterion for F:\\
\( F\circ S \circ G^{-1} \sqsubseteq R \;\;\; \Longleftrightarrow
\;\;\;  F \circ S \sqsubseteq R \rhd G \;\;\; \Longleftrightarrow
\;\;\;  F \sqsubseteq (R \rhd G)\rhd S^{-1} \) \\[0.8em]  
(2): Criterion for G:\\
\( F\circ S \circ G^{-1} \sqsubseteq R \;\;\;   \Longleftrightarrow
\;\;\;  (G \circ (F \circ S)^{-1})^{-1} \sqsubseteq  R \;\;\; 
\Longleftrightarrow \;\;\; ( G \circ S^{-1} F^{-1})^{-1} \sqsubseteq R     
 \;\;\;   \Longleftrightarrow \;\;\; G\circ S^{-1} \circ F^{-1}
 \sqsubseteq R^{-1}  \;\;\;   \Longleftrightarrow 
\;\;\;  G\circ S^{-1} \sqsubseteq R{-1} \rhd F  \;\;\;   \Longleftrightarrow
\;\;\; G \sqsubseteq (R^{-1} \rhd F)\rhd S    \)\\

\subsection{Translation into Fuzzy Monoidal Logics}
 We have seen that the lattice semantics can be translated easily into
 first order logic formulas of fuzzy monoidal logics. Let us look at
 the translation of some important properties of residuated lattices
 (defined above by Def. 3) into the 1st order logic formulas. These
 are listed in  
 Table 5. 

\begin{table}[t] 
\textbf{ 
\caption{\bf Lattice Semantics of First Order Logic Formulas of
  Monoidal fuzzy Logics}}\vspace{1em}  
\begin{align*} 
\mathbf{Lattice\;\; Semantics} &  & &
\mathbf{\;\;\;\;1st\;\;Order\;\;Logic\;\;Formulas}  
\end{align*}
\begin{align} 
 (x\otimes \bigvee_{i\in I} y_{i}) = \bigvee_{i\in I}& (x\otimes y)&
(x\& \exists y) = \exists y (x\& y) 
\end{align}
\begin{align}
 (x \Rightarrow \bigwedge_{i\in I} y_{i})=  \bigwedge_{i\in I}&
 (x\Rightarrow y_{i})&
(x \rightarrow \forall y) = \forall  (x\rightarrow y)
\end{align}
\begin{align}
\bigvee_{i\in I} x_{i}  \Rightarrow y  = \bigwedge_{i\in I} (x_{i}
 \Rightarrow y) &
 &(\exists x  \rightarrow y)  = \forall (x \rightarrow y) 
\end{align}
\begin{align}
 (x \otimes \bigwedge_{i\in I} y_{i})\leq  \bigwedge_{i\in I}&
 (x\otimes y_{i})&
(x \& \forall y) \rightarrow \forall  (x\& y)
\end{align} 
\begin{align}
\bigvee_{i\in I}( x  \Rightarrow y_{i})  \leq x 
 \Rightarrow \bigvee_{i\in I} y_{i}) &
 &\exists (x  \rightarrow y)  \rightarrow ( x \rightarrow \exists y) 
\end{align}
\begin{align}
\bigvee_{i\in I}( x_{i}  \Rightarrow y)  \leq \bigwedge_{i\in I} x_{i}
 \Rightarrow y &
 &\exists (x  \rightarrow y)  \rightarrow (\forall x \rightarrow y) 
\end{align} 
\label{Tab5}
\end{table}

More information about Monoidal Fuzzy Algebras and Monoidal Fuzzy
Logics can be found in \cite{hohle.crm},\cite{belohl.frs}. Monoidal logics
were originated by Ulrich H{\"{o}}hle. 

%Finally, we shall list the translation of the definitions of
%standard product $\circ$ and BK-products into first order quantified
%formulas: 

\section{Conclusion}
We have focused this paper towards examining some notions 
and technical features of non-associative compositions of mathematical 
relations that are fundamental in the logic of fuzzy relations and
also useful in applications. Many-valued logic based (fuzzy)
extensions of relations can  contribute on the theoretical 
side, by utilizing the elegant algebraic structure of relational systems. 
On the theoretical side,  fuzzy relations are extensions of standard  
non-fuzzy (crisp) relations. By replacing the usual Boolean algebra  
by many-valued logic algebras, one obtains extensions that contain  
the classical relational theory as a special case.  There is a whole
spectrum of systems covered by the structures presented in this
paper. 

What we have done is one coherent theory which still leaves  
the logician leeway to choose a specific base many valued  
logic algebra. One consistent algebraic meta-system which leads to  
formulas of great variety and allows for any number of specializations.  
As the general algebraic structure of relations has only   
minimal ontological commitment, this leaves also the engineer,
mathematician or scientist  with  choice (leeway) within which
different  fuzzy logics and ontologies can find elbow room.  
A number of different attitudes and needs can find space under this umbrella.

\end{document}